\documentclass[12pt]{article}
\usepackage{amsmath,amssymb}

\allowdisplaybreaks[4]
\newtheorem {Lemma}{Lemma}
\newtheorem {Theorem} {Theorem}

\newtheorem {Corollary}{Corollary}
\newenvironment {Proof} {\noindent {\it Proof }}{\hspace*{\fill}$\square$\par\vspace{3mm}}

\def\r{\mbox{}^r\mskip-4mu D'}

\date{}

\title{Reverse degree distance of unicyclic graphs}

\author{Zhibin Du, Bo Zhou\footnote{Corresponding author.}\\
Department of Mathematics, South China Normal University,\\
Guangzhou 510631,  China\\
e-mail: {\tt zhibindu@126.com; zhoubo@scnu.edu.cn}
}

\begin{document}

\maketitle

\begin{abstract}

The reverse degree distance is a connected graph invariant closely related to the degree distance proposed in mathematical chemistry. We
determine the unicyclic graphs of given girth, number of pendant
vertices and maximum degree, respectively, with maximum reverse
degree distances. 
\end{abstract}

\section{Introduction}

Let $G$  be a simple connected graph with vertex set $V(G)$.  For
$u,v\in V(G)$, let $d_G(u,v)$ be the distance between the vertices
$u$ and $v$ in $G$ and let $D_G(u)=\sum\limits_{v\in V(G)}d_G(u,v)$.
For $u\in V(G)$, let $d_G(u)$ be the degree of $u$ in $G$. The
degree distance of $G$ is defined as \cite{DK,ERA,Gu}
\[
D'(G)=\sum_{u\in V(G)}d_G(u)D_G(u).
\]
It is a useful molecular descriptor \cite{TCM}.
Earlier as noted in \cite{MNT,SSS}, this graph invariant appeared to
be part of the molecular topological index (or Schultz index)
\cite{Sch}, which may be expressed as $D'(G)+\sum\limits_{u\in
V(G)}d_G(u)^2$, see \cite{Gu,KG,MSKT,Zh}, where the latter part
$\sum\limits_{u\in V(G)}d_G(u)^2$ is known as the first Zagreb index
\cite{GRTW,GT,KKMT}. Thus, the degree distance is also called the
true Schultz index in chemical literature \cite{DGT}.


I. Tomescu \cite{To1} showed that the star is the unique graph with
minimum degree distance in the class of connected graphs with $n$
vertices. Further work on the minimum degree distance (especially for unicyclic and bicyclic graphs) may be found
in  A.I. Tomescu \cite{To2}, I. Tomescu \cite{To3} and Bucicovschi
and Cioab\v{a} \cite{BC}. Dankelmann et al.~\cite{DGMS}
gave asymptotically sharp upper bounds for the degree distance. Among others, the authors \cite{DZ2} studied the ordering of unicyclic graphs with large degree distances, and bicyclic graphs were also considered in \cite{DZ3}.



Recall that the Wiener index \cite{Wi} of the graph $G$ is defined
as
\[
W(G)=\frac{1}{2}\sum\limits_{u\in V(G)}D_G(u). \] Gutman \cite{Gu}
showed that if $G$ is a tree with $n$ vertices, then
\[D'(G)=4W(G)-n(n-1).
\] Thus there is no need to study the degree distance for trees
because this is equivalent to the study of the Wiener index, see,
e.g., \cite{DGKZ,DEG}.

The reverse degree distance of the graph $G$ is defined as \cite{ZT}
\[
\r (G)=2(n-1)md-D'(G)
\]
where $n$, $m$ and $d$ are the number of vertices, the number of
edges and the diameter of $G$, respectively.  Recall that, earlier,
Balaban et al.~\cite{BMIB} introduced the concept of reverse
Wiener index, which  is defined to be $\frac{n(n-1)d}{2}-W(G)$. If
$G$ is a tree, then from the result of Gutman \cite{Gu} mentioned
above,
\[
\r(G)=4\left[\frac{(n-1)^2d}{2}-W(G)\right]+n(n-1).
\]
Some properties of the reverse degree distance, especially for
trees, have been given in \cite{ZT}.  There are two reasons for the
study of this graph invariant.  One is that the reverse degree
distance itself is a topological index satisfying the basic
requirement to be a branching index and with potential for
application in chemistry \cite{ZT}. The other is the study the
reverse degree distance is actually the study the degree distance,
which is important in both mathematical chemistry and in discrete
mathematics.

In this paper, we determine the graphs with maximum reverse degrees
distance in the class of unicyclic graphs (connected graphs with a
unique cycle) with given girth (cycle length), number of pendant
vertices (vertices of degree one), and maximum degree, respectively.
Additionally, we also determine the graphs with minimum degree
distance in the class of unicyclic graphs with given number of
vertices, girth and diameter.

\section{Preliminaries}

Let $G$ be a graph of the form in Fig. 1, where $M$ and $N$ are
vertex-disjoint connected graphs, $T$ is a tree on $k \geq 2$
vertices such that $M$ and $T$ have only one common vertex $u$, and
$T$ and $N$ have only one common vertex $v$. Let $G^*$ be the graph
obtained from $M$ and $N$ by identifying vertices $u$ and $v$ which
is denoted by $u$,
and attaching $k-1$ pendant vertices to $u$.

\vspace{3mm}
\newpage

\begin{center}
\begin{picture}(30,-50)
\put(-80,-50){\begin{picture}(-60,-50)
\put(-20,21){\oval(39,25)}\put(0,20){\circle*{4.5}}\put(20,20){\circle{40}}
\put(40,20){\circle*{4.5}} \put(60,20){\oval(39,25)}
\put(-24,17){$M$}\put(57,17){$N$}\put(17,17){T}\put(18,-20){$G$}\put(-3,3){$u$}\put(37,3){$v$}
\end{picture}}
\end{picture}
\end{center}

\begin{center}
\begin{picture}(30,-50)
\put(90,-20){\begin{picture}(-60,-50)
\put(-20,18){\oval(39,25)}\put(0,20){\circle*{4.5}}
\put(0,20){\line(2,3){15}}\put(0,20){\line(-2,3){15}}
\put(-16,45){\circle*{4.5}}\put(16,45){\circle*{4.5}}
\put(-7,45){\circle*{2}}\put(0,45){\circle*{2}}\put(7,45){\circle*{2}}
\put(20,18){\oval(39,25)}
\put(-24,15){$M$}\put(18,15){$N$}\put(-4,-20){$G^*$}

\put(-13,59){${k-1}$} \put(-17,50){$\overbrace{~~~~~~~~~}$}

\put(-3,3){$u$}
\end{picture}}
\end{picture}
\end{center}

\vspace{10mm}
\begin{center}
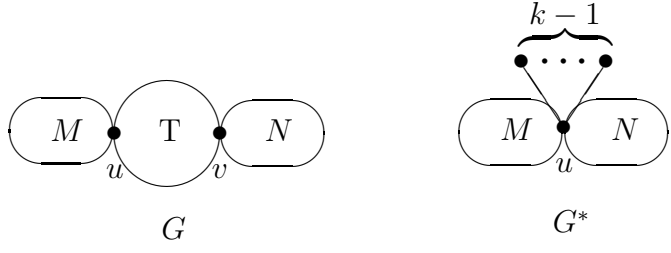
Fig.~1. The graphs $G$ and $G^*$.
\end{center}

\begin{Lemma}\label{lc}
Let $G$ and $G^*$ be the two graphs in Fig. 1.

\begin{enumerate}

\item[$(i)$]
If $V(N)=\{v\}$ and $G\not\cong G^*$, then $D'(G)> D'(G^*)$.

\item[$(ii)$]
If $|V(M)|,~|V(N)|\geq 3$, then $D'(G)>D'(G^*)$.
\end{enumerate}
\end{Lemma}

\begin{Proof}
For vertex--disjoint connected graphs  $Q_1$ and $Q_2$ with
$|V(Q_1)|$, $|V(Q_2)|\ge 2$, and $s\in V(Q_1)$, $t\in V(Q_2)$,  let
$H$ be the graph obtained from $Q_1$ and $Q_2$ by joining $s$ and
$t$ by an edge, and $H_1$ the graph obtained by identifying vertices
$s$ and $t$ which is denoted by $s$, and attaching a pendant vertex
$w$ to $s$.

Let $d_x=d_H(x)$ for $x\in V(H)$. It is easily seen that
\begin{eqnarray*}
&&(d_s+d_t-1)D_{H_1}(s)+1\cdot D_{H_1}(w)-d_sD_H(s)-d_tD_H(t)\\
&=&d_s[D_{H_1}(s)-D_H(s)]+d_t[D_{H_1}(s)-D_H(t)]+[D_{H_1}(w)-D_{H_1}(s)]\\
&=&-d_s(|V(Q_2)|-1)-d_t(|V(Q_1)|-1)+(|V(Q_1)|+|V(Q_2)|-2)\\
&=&-(d_s-1)(|V(Q_2)|-1)-(d_t-1)(|V(Q_1)|-1).
\end{eqnarray*}
Then
\begin{eqnarray*}
&&D'(H_1)-D'(H)\\
&=&-(|V(Q_2)|-1)\sum_{x\in V(Q_1)\setminus
\{s\}}d_x-(|V(Q_1)|-1)\sum_{x\in
V(Q_2)\setminus \{t\}}d_x\\
&&+(d_s+d_t-1)D_{H_1}(s)+1\cdot D_{H_1}(w)-d_sD_H(s)-d_tD_H(t)\\
&=&-(|V(Q_2)|-1)\sum_{x\in V(Q_1)\setminus
\{s\}}d_x-(|V(Q_1)|-1)\sum_{x\in
V(Q_2)\setminus \{t\}}d_x\\
&&-(d_s-1)(|V(Q_2)|-1)-(d_t-1)(|V(Q_1)|-1)<0,
\end{eqnarray*}
and thus $D'(H_1)<D'(H)$.

Now (i) and (ii) follow by  applying to $G$ the transformation from
$H$ to $H_1$ repeatedly.
\end{Proof}

\begin{Lemma}\label{l2}
Let $G_0$ be a connected graph with at least three vertices and let
$u$ and $v$ be two distinct vertices of $G_0$. Let $G_{s,t}$ be the
graph obtained from  $G_0$ by attaching $s$ and $t$ pendant vertices
to $u$ and $v$, respectively. If $s, t\ge 1$, then
$D'(G_{s,t})>\min\{D'(G_{s+t,0}), D'(G_{0,s+t})\}$.
\end{Lemma}

\begin{Proof}
Let $d_x=d_{G_0}(x)$ and $d(x,y)=d_{G_0}(x,y)$ for  $x,y\in V(G_0)$.
It is easily seen that
\begin{eqnarray*}
&&\left[(d_u+s+t)D_{G_{s+t,0}}(u)-(d_u+s)D_{G_{s,t}}(u)\right]+\left[d_vD_{G_{s+t,0}}(v)-(d_v+t)D_{G_{s,t}}(v)\right]\\
&=&(d_u+s)[D_{G_{s+t,0}}(u)-D_{G_{s,t}}(u)]+t[D_{G_{s+t,0}}(u)-D_{G_{s,t}}(v)]\\
&&+d_v[D_{G_{s+t,0}}(v)-D_{G_{s,t}}(v)]\\
&=&-t\cdot d(u,v)\cdot(d_u+s)+t\left[-sd(u,v)+\sum_{x\in
V(G_0)\setminus
\{u,v\}}(d(x,u)-d(x,v))\right]\\
&&+t\cdot d(u,v)\cdot d_v\\
&=&t\left[(d_v-d_u-2s)d(u,v)+\sum_{x\in
V(G_0)\setminus \{u,v\}}(d(x,u)-d(x,v))\right]\\
\end{eqnarray*}
and thus
\begin{eqnarray*}
&&D'(G_{s+t,0})-D'(G_{s,t})\\
&=&t\sum_{x\in V(G_0)\setminus \{u,v\}}d_x\left(d(x,u)-d(x,v)\right)\\
&&-std(u,v)+t\bigg[-sd(u,v)+\sum_{x\in
V(G_0)\setminus \{u,v\}}(d(x,u)-d(x,v))\bigg]\\
&&+(d_u+s+t)D_{G_{s+t,0}}(u)-(d_u+s)D_{G_{s,t}}(u)\\
&&+d_vD_{G_{s+t,0}}(v)-(d_v+t)D_{G_{s,t}}(v)\\
&=&t\sum_{x\in
V(G_0)\setminus \{u,v\}}(d_x+1)(d(x,u)-d(x,v))-2std(u,v)\\
&&+t\left[(d_v-d_u-2s)d(u,v)+\sum_{x\in
V(G_0)\setminus \{u,v\}}(d(x,u)-d(x,v))\right]\\
&=&t\bigg[(d_v-d_u-4s)d(u,v)+\sum_{x\in V(G_0)\setminus
\{u,v\}}(d_x+2)(d(x,u)-d(x,v))\bigg].
\end{eqnarray*}
Similarly, we have
\begin{eqnarray*}
&&D'(G_{0,s+t})-D'(G_{s,t})\\
&=&s\bigg[(d_u-d_v-4t)d(u,v)+\sum_{x\in V(G_0)\setminus
\{u,v\}}(d_x+2)(d(x,v)-d(x,u))\bigg].
\end{eqnarray*}
If $D'(G_{s+t,0})\ge D'(G_{s,t})$, then
\begin{eqnarray*}
\sum_{x\in V(G_0)\setminus \{u,v\}}(d_x+2)(d(x,v)-d(x,u))&\le&
(d_v-d_u-4s)d(u,v)
\end{eqnarray*}
and thus,
\begin{eqnarray*}
D'(G_{0,s+t})-D'(G_{s,t})&\le&s\bigg[(d_u-d_v-4t)d(u,v)+(d_v-d_u-4s)d(u,v)\bigg]\\
&=&-4s(s+t)d(u,v)<0.
\end{eqnarray*}
The result follows.
\end{Proof}

Let $G$ and $H$ be connected graphs. Let $V_1(G)=\{x\in
V(G):d_G(x)=2\}$ and $V_2(G)=V(G)\setminus V_1(G)$.  Let
$d_x=d_G(x)$ for $x\in V(G)$, and $d^*_x=d_H(x)$ for $x\in V(H)$.
Then
\begin{eqnarray*}
&&D'(H)-D'(G)\\
&=&2\sum_{x\in V_1(H)}D_H(x)+\sum_{x\in
V_2(H)}d^*_xD_H(x)-2\sum_{x\in V_1(G)}D_{G}(x)-\sum_{x\in
V_2(G)}d_xD_{G}(x)\\
&=&4[W(H)-W(G)]+\sum_{x\in V_2(H)}(d^*_x-2)D_H(x)-\sum_{x\in
V_2(G)}(d_x-2)D_{G}(x).
\end{eqnarray*}

Let $G^*$ be the unicyclic graph obtained from the cycle
$C_m=v_0v_1\dots v_{m-1}v_0$ by attaching a path $P_a$ and a path
$P_b$ to $v_i$ and $v_j$, respectively, where $i\ne j$, $a\ge 1$ and
$b\ge 2$. Label the vertices of the path $P_b$  attached to $v_j$ as
$u_1,\dots,u_b$ consecutively, where $u_1$ is adjacent to $v_j$ in
$G^*$.

For integer $h\ge 1$,  let $G_{u_t,h}^{(1)}$ be the graph obtained
from $G^*$ by attaching $h$ pendant vertices to $u_t$, where $1\le
t\le b-1$,  and $G_{v_t,h}^{(2)}$ the graph obtained from $G^*$ by
attaching $h$ pendant vertices to $v_t$, where $0\le t\le m-1$.

\begin{Lemma}\label{l5}
Let $n_1=a+m-1$ and $n_2=b-t$ in $G_{u_t,h}^{(1)}$. Then
$D'\left(G^{(2)}_{v_j,h}\right)-D'\left(G^{(1)}_{u_t,h}\right)=2ht[2(n_2-n_1)-1]$.
\end{Lemma}

\begin{Proof}
Let $v$ be a pendant vertex attached to $v_j$ in $G^{(2)}_{v_j,h}$
(resp. $u_t$ in $G^{(1)}_{u_t,h}$). Let $G_1=G^{(2)}_{v_j,h}$ and
$G_2=G^{(1)}_{u_t,h}$. Then
\begin{eqnarray*}
&&D'\left(G^{(2)}_{v_j,h}\right)-D'\left(G^{(1)}_{u_t,h}\right)\\
&=&4[W(G_1)-W(G_2)]-[D_{G_1}(u_b)-D_{G_2}(u_b)]\\
&&-h[D_{G_1}(v)-D_{G_2}(v)]+(h+1)D_{G_1}(v_j)-hD_{G_2}(u_t)-D_{G_2}(v_j)\\
&=&4ht(n_2-n_1)-ht-ht(n_2-n_1)+ht(n_2-n_1)-ht\\
&=&2ht[2(n_2-n_1)-1],
\end{eqnarray*}
as desired.
\end{Proof}

By similar arguments as in Lemma \ref{l5},  we have

\begin{Lemma}\label{l6}
Let $c=d_G(v_i,v_j)$, $t_1=d_G(v_i,v_t)$ and $t_2=d_G(v_j,v_t)$,
where  $G=G^{(2)}_{v_t,h}$. If $t\ne i$, then
$D'\left(G^{(2)}_{v_i,h}\right)-D'\left(G^{(2)}_{v_t,h}\right)=4h[b(c-t_2)-at_1]$.
\end{Lemma}

As usual, $G-V_1$ means the graph formed from the graph $G$ by
deleting the vertices of $V_1\subset V(G)$ and edges incident with
these vertices, while $G-E_1$ means the graph formed from $G$ by
deleting edges of $E_1\subseteq E(G)$.

\section{Minimum Degree distance of unicyclic graphs with given girth and diameter}

In this section we determine the unicyclic graphs with minimum
degree distance when the number of vertices, girth and diameter are
given.

Let $n$, $m$ and $d$ be integers with  $3\le m\le n-1$ and $2\le d
\le n-\lfloor\frac{m+1}{2}\rfloor$. Let $P_s$ be a path on $s$
vertices. For $a\ge b\ge 0$ and $a\ge 1$, let $U^k_{n,m,d}(a,b)$ be
the unicyclic graph obtained from the cycle $C_m=v_0v_1\dots
v_{m-1}v_0$ by attaching a path $P_a$ to $v_0$ and a path $P_b$ to
$v_{\lfloor\frac{m}{2}\rfloor}$ respectively (if $b=0$, then by
attaching only a path $P_a$ to $v_0$), where
$a+b=d-\lfloor\frac{m}{2}\rfloor$, and attaching
$n-d-\lfloor\frac{m+1}{2}\rfloor$ pendant vertices to $v_k$, where
$0\le k\le \lfloor\frac{m}{4}\rfloor$. Let
$U_{n,m,d}(a,b)=U^0_{n,m,d}(a,b)$.

For $U_{n,m,d}(a,b)$, let $u_0$ be the pendant vertex on the path
attached to $v_0$, let $u_1$ be the pendant vertex on the path
attached to $v_{\lfloor\frac{m}{2}\rfloor}$ if $b\ge 1$, and
$u_1=v_{\lfloor\frac{m}{2}\rfloor}$ if $b=0$, let $u$ be any of the
pendant vertices attached to $v_0$.


Let $\alpha=\alpha(n,m,d)=
\frac{\left(n-d-\left\lfloor\frac{m+1}{2}\right\rfloor\right)\left\lfloor\frac{m}{2}\right\rfloor}{n-d-\frac{1}{2}}$.
Let $\gamma$ and $\theta$ be integers such that
$\gamma+\theta=d-\left\lfloor\frac{m}{2}\right\rfloor$ and
$\gamma-\theta$ is an integer as large as possible but no more than
$\alpha+1$. Let $U_{n,m,d}=U_{n,m,d}(\gamma,\theta)=
U_{n,m,d}^0(\gamma,\theta)$.

\begin{Lemma} \label{add} Let  $n$, $m$ and $d$ be fixed integers with $3\le m\le n-2$ and $3\le d
\le n-\lfloor\frac{m+1}{2}\rfloor$. Then
$D'\left(U_{n,m,d}(a,b)\right)$ with $a\ge b$ and
$a+b=d-\left\lfloor\frac{m}{2}\right\rfloor$ is minimum if and only
if $(a, b)=(\gamma, \theta)$, $(\gamma-1,\theta+1)$ if $\alpha\ge 1$
is an integer with different parity as
$d-\left\lfloor\frac{m}{2}\right\rfloor$, and $(a, b)=(\gamma,
\theta)$ otherwise.
\end{Lemma}

\begin{Proof}
Let $h=n-d-\lfloor\frac{m+1}{2}\rfloor$. Let $w$ be the neighbor of
$u_0$ in $U_{n,m,d}(a,b)$. Note that for $a-b\ge2$,
$U_{n,m,d}(a-1,b+1)\cong U_{n,m,d}(a,b)-\{u_0\}+\{u_0u_1\}$. Let
$G_1=U_{n,m,d}(a-1,b+1)$ and $G_2=U_{n,m,d}(a,b)$. If $a\ge b\ge 1$,
then
\begin{eqnarray*}
&&D'(U_{n,m,d}(a-1,b+1))-D'(U_{n,m,d}(a,b))\\
&=&4[W(G_1)-W(G_2)]-h[D_{G_1}(u)-D_{G_2}(u)]-[D_{G_1}(u_0)-D_{G_2}(u_0)]\\
&&+[D_{G_1}(v_{\lfloor\frac{m}{2}\rfloor})-D_{G_2}(v_{\lfloor\frac{m}{2}\rfloor})]+(h+1)[D_{G_1}(v_0)-D_{G_2}(v_0)]\\
&&-D_{G_1}(w)+D_{G_2}(u_1)\\
&=&4\left[(1-a+b)\left(h+\left\lfloor\frac{m-1}{2}\right\rfloor+\frac{1}{2}\right)
+h\left\lfloor\frac{m}{2}\right\rfloor\right],
\end{eqnarray*}
and if $a=d-\lfloor\frac{m}{2}\rfloor$ and $b=0$, then
\begin{eqnarray*}
&&D'(U_{n,m,d}(a-1,b+1))-D'(U_{n,m,d}(a,b))\\
&=&4[W(G_1)-W(G_2)]-h[D_{G_1}(u)-D_{G_2}(u)]-[D_{G_1}(u_0)-D_{G_2}(u_0)]\\
&&+(h+1)[D_{G_1}(v_0)-D_{G_2}(v_0)]+D_{G_1}(v_{\lfloor\frac{m}{2}\rfloor})-D_{G_1}(w)\\
&=&4\left[(1-a)\left(h+\left\lfloor\frac{m-1}{2}\right\rfloor+\frac{1}{2}\right)
+h\left\lfloor\frac{m}{2}\right\rfloor\right].
\end{eqnarray*}
Thus $D'(U_{n,m,d}(a-1,b+1))\ge D'(U_{n,m,d}(a,b))$ if and only if
$a-b\le \alpha+1$, and $D'(U_{n,m,d}(a-1,b+1))=D'(U_{n,m,d}(a,b))$
if and only if $a-b=\alpha+1$. Thus $D'(U_{n,m,d}(a,b))$ is minimum
if and only if $a-b$ is as large as possible with $a-b\le \alpha+1$.
Note that $a-b=\alpha +1$ if and only if $\alpha\ge 1$ is an integer
with different parity as $d-\left\lfloor\frac{m}{2}\right\rfloor$.
The result follows.
\end{Proof}

Let $\mathbb{U}(n,m,d)$ be the set of unicyclic graphs with $n$
vertices, girth $m$ and diameter $d$, where $2\le d\le
n-\lfloor\frac{m+1}{2}\rfloor$ and $3\le m\le n-2$. If $G\in
\mathbb{U}(n,m,2)$,  then $m=3$ and $G=U_{n,3,2}(1,0)$.

Let $G$ be a unicyclic graph with $n$ vertices and let
$C_m=v_0v_1\dots v_{m-1}v_0$ be its unique cycle. Then $G-E(C_m)$
consists of $m$ trees $T_0, T_1, \dots, T_{m-1}$, where $v_i\in
V(T_i)$ for $i=0,1,\dots, m-1$. If the degree of $v_i$ is at least
three, then the components of $T_i-v_i$ are called the branches of
$G$ at $v_i$, each containing a neighbor of $v_i$ in $T_i$.

\begin{Lemma} \label{newadd}
Let $n$, $m$ and $d$ be integers with $n\ge6$, $3\le m\le n-2$ and
$3\le d \le n-\lfloor\frac{m+1}{2}\rfloor$, and let
$\beta=\frac{1}{2}(d-\lfloor\frac{m}{2}\rfloor)$. If $G$ is a graph
with minimum degree distance in $\mathbb{U}(n,m,d)$, then
$G=U_{n,m,d}(a,b)=U_{n,m,d}^0(a,b)$ with $a \ge b$ or
$G=U_{n,m,d}^k(\beta,\beta)$ for
$k=1,\dots,\lfloor\frac{m}{4}\rfloor$.
\end{Lemma}

\begin{Proof}
Let $C_m=v_0v_1\dots v_{m-1}v_0$ be the unique cycle of $G$, and let
$P(G)=u_0u_1\dots u_d$ be a diametrical path of $G$.  Let
$d(x,y)=d_G(x,y)$ for $x,y\in V(G)$.

Suppose that $P(G)$ has no common vertices with $C_m$. Let $u_s$ and
$v_t$ be the vertices such that $d(u_s, v_t)=\min\{d(u,v): u\in
V(P(G)), v\in V(C_m)\}$. Using Lemma \ref{lc} (ii) by setting
$u=u_s$, $v=v_t$, $M$ to be the subgraph of $G$ consisting of the
path $P(G)$ and trees attached to $u_i$ for all $1\le i\le d-1$ and
$i\ne s$, $N$ to be the subgraph of $G$ by deleting all branches at
$v_t$, we obtain a graph $G^*$ for which $P(G^*)$ $(=P(G))$ and the
cycle $C_m$ have exactly one common vertex and $D'(G^*) < D'(G)$, a
contradiction. Thus, $P(G)$ and $C_m$ have at least one common
vertex. We may choose $P(G)$ such that $P(G)$ and the cycle $C_m$
have vertices in common as many as possible and $u_0$ is a pendant
vertex.

Let $v_i=u_a$ (resp. $v_j=u_l$) be the first (resp. last) common
vertex of $P(G)$ and $C_m$, where $0<a\le l\le d$. By Lemma \ref{lc}
(i), all vertices outside $C_m$ except those in $T_i$ and $T_j$ are
pendant vertices attached to vertices that are nearest to them in
$C_m$, all vertices in $T_i$ and $T_j$ except those in $P(G)$ are
pendant vertices attached to vertices that are nearest to them in
$P(G)$.

Suppose that $P(G)$ and $C_m$ have only one common vertex, i.e.,
$i=j$, $a=l$ and $l<d$. By the choice of $P(G)$, we have $a\ge2$. By
Lemma \ref{l2}, all pendant vertices in $G$ except $u_0$ and $u_d$
are actually attached to some vertex, say $s$, of $G$.

Suppose that $s\in\{v_0,v_1, \dots, v_{m-1}\}\setminus\{v_i\}$, say
$s=v_q$. Let $N_q=\{v_{q_1},\dots,v_{q_t}\}$ be the set of pendant
vertices attached to $v_q$. For
$H=G-\{v_qv_{q_1},\dots,v_qv_{q_t}\}+\{v_iv_{q_1},\dots,v_iv_{q_t}\}\in\mathbb{U}(n,m,d)$,
we have
\begin{eqnarray*}
&&D'(H)-D'(G)\\
&=&4[W(H)-W(G)]-[D_H(u_0)-D_G(u_0)]-[D_H(u_d)-D_G(u_d)]\\
&&-t[D_H(v_{q_1})-D_G(v_{q_1})]+(t+2)D_H(v_i)-tD_G(v_q)-2D_G(v_i)\\
&=&-4dt\cdot d(v_q,v_i)+t\cdot d(v_q,v_i)+t\cdot d(v_q,v_i)\\
&&+dt\cdot d(v_q,v_i)-dt\cdot d(v_q,v_i)-2t\cdot d(v_q,v_i)\\
&=&-4dt\cdot d(v_q,v_i),
\end{eqnarray*}
and then $D'(H)<D'(G)$, a contradiction. Thus, $s\in \{u_1,
u_2,\dots, u_{d-1}\}$. Suppose without loss of generality that
$s\in\{u_a, u_{a+1}\dots,u_{d-1}\}$. For
$H^{*}=G-\{u_{a-2}u_{a-1}\}+\{v_{i-1}u_{a-2}\}\in\mathbb{U}(n,m,d)$,
the path $P(H^{*})=u_0\dots u_{a-2}v_{i-1}v_iu_{a+1}\dots u_d$ has
more than one common vertex with the cycle $C_m$ and the same length
as $P(G)$, and we have
\begin{eqnarray*}
D'(H^{*})-D'(G)&=&4[W(H^{*})-W(G)]-[D_{H^{*}}(u_0)-D_G(u_0)]\\
&&+D_{H^{*}}(v_{i-1})-D_{H^{*}}(u_{a-1})\\
&=&-4(m-2)(a-1)+(m-2)-2a-m+4\\
&=&-2(a-1)(2m-3)<0,
\end{eqnarray*}
and then $D'(H^{*})<D'(G)$, a contradiction. Thus, $P(G)$ and $C_m$
have at least two common vertices, i.e., $a<l$.

By Lemma \ref{l2}, all pendant vertices in $G$ except $u_0$ and
$u_d$ are actually attached to some vertex, say $x$, in $G$. Thus
$x$ has exactly $h=n-m-a-(d-l)$ pendant neighbors outside $P(G)$.
Let $b=d-l$. Assume that $a\ge b$.

Suppose that $l<d$ and $x\in\{u_{l+1},\dots,u_{d-1}\}$, say $x=u_q$,
where $l<q\leq d-1$. Let $u_{q_1},u_{q_2},\dots,u_{q_h}$ be the
pendant neighbors of $u_q$ outside $P(G)$. For
$G_1=G-\{u_qu_{q_1},u_qu_{q_2},\dots,u_qu_{q_h}\}+\{u_lu_{q_1},u_lu_{q_2},\dots,u_lu_{q_h}\}\in\mathbb{U}(n,m,d)$,
using Lemma \ref{l5} by setting $t=q-l$, $n_1=a+m-1$ and $n_2=b-t$,
and noting that $n_1>n_2$ since $a\ge b$, we have
\[
D'(G_1)-D'(G)=2h(q-l)[2(n_2-n_1)-1]<0,
\]
and then $D'(G_1)<D'(G)$, a contradiction. Thus,
$x\not\in\{u_{l+1},u_{l+2}, \dots,u_{d-1}\}$ if $l<d$. Moreover, if
$a=b$ then by similar arguments,
$x\not\in\{u_{1},u_2\dots,u_{a-1}\}$, and thus $x\in\{v_0,v_1,
\dots, v_{m-1}\}$.

\noindent {\bf Case 1.} $a>b$.

First we prove that $x\in \{u_1,u_2,\dots,u_a\}$. Suppose that this
is not true.  Then  $x=v_s$ for some $s$ with $0\le s\le m-1$ and
$s\ne i$. Let $N_s=\{v_{s_1},\dots,v_{s_h}\}$ be the set of pendant
vertices attached to $v_s$. Suppose that $d(v_i,v_j)=c$,
$d(v_i,v_s)=t_1$ and $d(v_j,v_s)=t_2$, then $c\le t_1+t_2$. Consider
$G_2=G-\{v_sv_{s_1},\dots, v_sv_{s_h}\}+\{v_iv_{s_1},\dots,
v_iv_{s_h}\}\in\mathbb{U}(n,m,d)$. Note that if $l=d$, then $b=0$.
By Lemma \ref{l6}, we have
\begin{eqnarray*}
D'(G_2)-D'(G)&=&4h[b(c-t_2)-at_1]\le 4h(bt_1-at_1)=4ht_1(b-a)<0,
\end{eqnarray*}
and then $D'(G_2)<D'(G)$, a contradiction. Thus, $x\in
\{u_1,u_2,\dots,u_a\}$, say $x=u_p$, where $1\leq p \leq a$.

Next we prove that $d(v_i, v_j)=\lfloor \frac{m}{2}\rfloor$. If
$l=d$, then this is obvious. Suppose that $l<d$ and $c=d(v_i,
v_j)<\lfloor\frac{m}{2}\rfloor$. Let $v$ be the neighbor of $v_j$ on
$C_m$ with $d(v_i, v)=c+1$ (If $\{v_i,v_{i+1},\dots,v_{j-1},v_j\}$
is the shortest path from $v_i$ to $v_j$, then $v=v_{j+1}$). By our
choice of $P(G)$, we have $b+c>\lfloor \frac{m}{2}\rfloor$, and then
$b>1$. Note that $c>0$. Consider
$G_3=G-\{v_ju_{l+1}\}+\{vu_{l+1}\}-\{u_{d-1}u_d\}+\{v_iu_d\}\in\mathbb{U}(n,m,d)$.
If $1\le p\le a-1$, then
\begin{eqnarray*}
&&D'(G_3)-D'(G)\\
&=&4[W(G_3)-W(G)]-[D_{G_3}(u_0)-D_G(u_0)]-[D_{G_3}(u_d)-D_G(u_d)]\\
&&+h[D_{G_3}(u_p)-D_G(u_p)]-h[D_{G_3}(u_{p_1})-D_G(u_{p_1})]\\
&&+2D_{G_3}(v_i)+D_{G_3}(v)-D_{G_3}(u_{d-1})-D_G(v_i)-D_G(v_j)\\
&=&-2(2m-3)(b-1)-4c(n-m-2b+1)<0.
\end{eqnarray*}
By similar calculation, $D'(G_3)-D'(G)=-2(2m-3)(b-1)-4c(n-m-2b+1)<0$
holds also if $p=a$.  It follows that in either case
$D'(G_3)<D'(G)$, a contradiction. Thus, $d(v_i, v_j)=\lfloor
\frac{m}{2}\rfloor$ and $h=n-d-\lfloor\frac{m+1}{2}\rfloor$.

Now we prove $p=a$. Suppose that $p\le a-1$. Let
$u_{p_1},u_{p_2},\dots,u_{p_h}$ be the pendant neighbors of $u_p$
outside $P(G)$. If $b+m>a$, then $p<b+m-1$, and for
$G_4=G-\{u_pu_{p_1},\dots,u_pu_{p_h}\}+\{u_au_{p_1},\dots,u_au_{p_h}\}\in
\mathbb{U}(n,m,d)$, using Lemma \ref{l5} by setting $t=a-p$,
$n_1=b+m-1$ and $n_2=p$, we have
\begin{eqnarray*}
D'(G_4)-D'(G)&=&2h(a-p)[2p-2(b+m-1)-1]<0,
\end{eqnarray*}
and thus  $D'(G_4)<D'(G)$, a contradiction. Suppose that $b+m\le a$.
Then $b-(h-1)<a$.  Consider $G_4=G-\{u_pu_{p_1},\dots,u_pu_{p_h}\}+
\{u_{p+1}u_{p_1},\dots,u_{p+1}u_{p_h}\}-\{u_0u_1\}+\{u_0u_d\}\in
\mathbb{U}(n,m,d)$. If $l<d$ and $1\le p\le a-2$, then
\begin{eqnarray*}
&&D'(G_4)-D'(G)\\
&=&4[W(G_4)-W(G)]+[D_{G_4}(u_a)-D_G(u_a)]+[D_{G_4}(u_l)-D_G(u_l)]\\
&&-[D_{G_4}(u_0)-D_G(u_0)]-h[D_{G_4}(u_{p_1})-D_G(u_{p_1})]\\
&&+hD_{G_4}(u_{p+1})-D_{G_4}(u_1)-hD_G(u_p)+D_G(u_d)\\
&=&2\left(2\left\lfloor\frac{m-1}{2}\right\rfloor+1\right)(b-a+1-h)<0.
\end{eqnarray*}
By similar calculation, the inequality $D'(G_4)-D'(G)<0$ holds also
if $l=d$ or $p=a-1$. Then, in any case, $D'(G_4)<D'(G)$, a
contradiction. Thus, $p=a$.

Now we have proved that $G=U_{n,m,d}(a,b)$, where $a>b$ and
$a+b=d-\lfloor\frac{m}{2}\rfloor$.

\noindent {\bf Case 2.} $a=b$.

Note that $x\in \{v_0, v_1,\dots, v_{m-1}\}$, say $x=v_s$. Assume
that $v_iv_{i+1}\dots v_{j-1}v_j$ is a shortest path from $v_i$ to
$v_j$.  Obviously, $m\ge2(j-i)$. If $m=2(j-i)$, then by symmetry, we
may assume that $i\le s\le j$. Suppose that $m> 2(j-i)$ and
$s\not\in\{i,i+1,\dots,j-1,j\}$. Let $N_s=\{v_{s_1},\dots,v_{s_h}\}$
be the set of pendant vertices attached to $v_s$. Let
$d(v_i,v_j)=c$, $d(v_i,v_s)=t_1$ and $d(v_j,v_s)=t_2$. Then
$c<t_1+t_2$.  For $G_5=G-\{v_sv_{s_1},\dots,
v_sv_{s_h}\}+\{v_iv_{s_1},\dots, v_iv_{s_h}\}\in\mathbb{U}(n,m,d)$,
by Lemma \ref{l6}, we have
\begin{eqnarray*}
D'(G_5)-D'(G)=4h[b(c-t_2)-at_1]=4ha(c-t_1-t_2)<0,
\end{eqnarray*}
and then $D'(G_5)<D'(G)$, a contradiction. Thus $i\le s\le j$.

Suppose that $c=d(v_i, v_j)<\lfloor\frac{m}{2}\rfloor$. Note that
$d(v_i,v_{j+1})=c+1$. By our choice of $P(G)$, we have $b+c>\lfloor
\frac{m}{2}\rfloor$, and then $b>1$. Consider
$G_6=G-\{v_ju_{l+1}\}+\{v_{j+1}u_{l+1}\}-\{u_{d-1}u_d\}+\{v_su_d\}\in\mathbb{U}(n,m,d)$.
If $i+1\le s\le j-1$, then
\begin{eqnarray*}
&&D'(G_6)-D'(G)\\
&=&4[W(G_6)-W(G)]-[D_{G_6}(u_0)-D_G(u_0)]-[D_{G_6}(u_d)-D_G(u_d)]\\
&&-h[D_{G_6}(v_{s_1})-D_G(v_{s_1})]+[D_{G_6}(v_i)-D_G(v_i)]\\
&&+(h+1)D_{G_6}(v_s)-D_{G_6}(u_{d-1})+D_{G_6}(v_{j+1})-hD_G(v_s)-D_G(v_j)\\
&=&-2(2m-3)(b-1)-4(j-s)(h+1)<0.
\end{eqnarray*}
By similar calculation, the inequality $D'(G_6)-D'(G),0$ holds also
if $s=i$ or $j$. In either case, we have $D'(G_6)<D'(G)$, a
contradiction. Thus $d(v_i,v_j)=\lfloor\frac{m}{2}\rfloor$ and
$h=n-d-\lfloor\frac{m+1}{2}\rfloor$. By Lemma \ref{l6}, we have
$U^k_{n,m,d}(\beta,\beta)$ for
$k=0,1,\dots,\lfloor\frac{m}{4}\rfloor$ have equal degree distance,
and thus $G= U^k_{n,m,d}(\beta,\beta)$ for
$k=0,1,\dots,\lfloor\frac{m}{4}\rfloor$.

By combining Cases 1 and 2, we have
$G=U_{n,m,d}(a,b)=U_{n,m,d}^0(a,b)$ with $a \ge b$ or
$G=U_{n,m,d}^k(\beta,\beta)$ for
$k=1,\dots,\lfloor\frac{m}{4}\rfloor$.
\end{Proof}

\begin{Theorem} \label{th1}
Let $n$, $m$ and $d$ be integers with $n\ge6$, $3\le m\le n-2$ and
$3\le d \le n-\lfloor\frac{m+1}{2}\rfloor$, and let
$\alpha=\frac{\left(n-d-\left\lfloor\frac{m+1}{2}\right\rfloor\right)\left\lfloor\frac{m}{2}\right\rfloor}{n-d-\frac{1}{2}}$,
$\beta=\frac{1}{2}(d-\lfloor\frac{m}{2}\rfloor)$.

\begin{enumerate}

\item[(i)]
If ~$0<\alpha<1$ and $d-\lfloor\frac{m}{2}\rfloor$ is even, then
$U^k_{n,m,d}(\beta,\beta)$ for
$k=0,1,\dots,\lfloor\frac{m}{4}\rfloor$ are the unique graphs in
$\mathbb{U}(n,m,d)$ with minimum degree distance.

\item[(ii)]
If ~$\alpha=1$ and $d-\lfloor\frac{m}{2}\rfloor$ is even, then
$U_{n,m,d}=U_{n,m,d}(\beta+1,\beta-1)$ and
$U^k_{n,m,d}(\beta,\beta)$ for
$k=0,1,\dots,\lfloor\frac{m}{4}\rfloor$ are the unique graphs in
$\mathbb{U}(n,m,d)$ with minimum degree distance.

\item[(iii)]
If ~$\alpha>1$ is an integer with different parity as
$d-\lfloor\frac{m}{2}\rfloor$, then $U_{n,m,d}$ and
$U_{n,m,d}(\gamma-1,\theta+1)$ are the unique graphs in
$\mathbb{U}(n,m,d)$ with minimum degree distance.

\item[(iv)]
If ~$\alpha=0$, or $0<\alpha\le1$ and $d-\lfloor\frac{m}{2}\rfloor$
is odd, or $\alpha>1$ is either not an integer or an integer with
the same parity as $d-\lfloor\frac{m}{2}\rfloor$, then $U_{n,m,d}$
is the unique graph in $\mathbb{U}(n,m,d)$ with minimum degree
distance.

\end{enumerate}
\end{Theorem}

\begin{Proof}
Suppose that $G$ is a graph with minimum degree distance in
$\mathbb{U}(n,m,d)$. By Lemma \ref{newadd}, we have
$G=U_{n,m,d}(a,b)=U_{n,m,d}^0(a,b)$ with $a \ge b$ or
$G=U_{n,m,d}^k(\beta,\beta)$ for
$k=1,\dots,\lfloor\frac{m}{4}\rfloor$. If $G=U_{n,m,d}^0(a,b)$, then
we have by Lemma \ref{add} that $G=U_{n,m,d}$ or
$U_{n,m,d}(\gamma-1,\theta+1)$ if $\alpha\ge 1$ is an integer with
different parity as $d-\left\lfloor\frac{m}{2}\right\rfloor$, and
$G=U_{n,m,d}$ otherwise.

If $d-\lfloor\frac{m}{2}\rfloor$ is odd, then  $G=U_{n,m,d}$,
$U_{n,m,d}(\gamma-1,\theta+1)$ if $\alpha\ge1$ is an even integer,
and $G=U_{n,m,d}$ otherwise.

Suppose that $d-\lfloor\frac{m}{2}\rfloor$ is even. Then either
$G=U^k_{n,m,d}(\beta,\beta)$ for
$k=1,2,\dots,\lfloor\frac{m}{4}\rfloor$, or
 $G=U_{n,m,d}$,
$U_{n,m,d}(\gamma-1,\theta+1)$ if $\alpha\ge1$ is an odd integer,
and $G=U_{n,m,d}$ otherwise.

If $\alpha=0$, then $G=U_{n,m,d}$.

If $0<\alpha<1$, then $G=U^k_{n,m,d}(\beta,\beta)$ for
$k=0,1,\dots,\lfloor\frac{m}{4}\rfloor$.

Suppose that $\alpha=1$. Then $G=U_{n,m,d}$, $U_{n,m,d}(\gamma-1,
\theta+1)$, or $G=U^k_{n,m,d}(\beta,\beta)$ for
$k=1,2,\dots,\lfloor\frac{m}{4}\rfloor$. Since $(\gamma-1,
\theta+1)=(\beta,\beta)$, we have $G=U_{n,m,d}$, or
$G=U^k_{n,m,d}(\beta,\beta)$ for
$k=0,1,\dots,\lfloor\frac{m}{4}\rfloor$.

Suppose that $\alpha>1$ is an odd integer. Then $G=U_{n,m,d}$,
$U_{n,m,d}(\gamma-1,\theta+1)$, or $G=U^k_{n,m,d}(\beta,\beta)$ for
$k=1,2,\dots,\lfloor\frac{m}{4}\rfloor$. Let
$h=n-d-\lfloor\frac{m+1}{2}\rfloor$. Then $\alpha=
\frac{h\lfloor\frac{m}{2}\rfloor}{h+\lfloor\frac{m-1}{2}\rfloor+\frac{1}{2}}$
and $h>0$. By the proof of Lemma \ref{add}, we have
\begin{eqnarray*}
&&D'(U^0_{n,m,d}(\beta,\beta))-D'(U_{n,m,d}(\gamma,\theta))\\
&=&2(\gamma-\theta)\left[\frac{1}{2}(\theta-\gamma)\left(h+\left\lfloor\frac{m-1}{2}\right\rfloor+\frac{1}{2}\right)+h\left\lfloor\frac{m}{2}\right\rfloor\right]\\
&>&2\left[-\frac{1}{2}(\alpha+1)\left(h+\left\lfloor\frac{m-1}{2}\right\rfloor+\frac{1}{2}\right)+h\left\lfloor\frac{m}{2}\right\rfloor\right]\\
&=&\frac{h(\alpha-1)}{\alpha}\left\lfloor\frac{m}{2}\right\rfloor>0,
\end{eqnarray*}
and then $D'(U^0_{n,m,d}(\beta,\beta))>D'(U_{n,m,d}(\gamma,\theta))=
D'\left(U_{n,m,d}(\gamma-1,\theta+1)\right)$. Thus $G=U_{n,m,d}$,
$U_{n,m,d}(\gamma-1,\theta+1)$.

Suppose that $\alpha>1$ is not an odd integer. Then $G=U_{n,m,d}$,
or $G=U^k_{n,m,d}(\beta,\beta)$ for
$k=1,2\dots,\lfloor\frac{m}{4}\rfloor$. By similar arguments as
above, we have
$D'\left(U^0_{n,m,d}(\beta,\beta)\right)>D'\left(U_{n,m,d}\right)$.
Thus $G=U_{n,m,d}$.
\end{Proof}

\begin{Corollary}
\label{Cor} Let $G\in \mathbb{U}(n,m,d)$  with $n\ge6$, $3\le m\le
n-2$ and $3\le d \le n-\lfloor\frac{m+1}{2}\rfloor$.  Then $D'(G)\ge
D'\left(U_{n,m,d}\right)$.
\end{Corollary}

\section{Reverse degree distances of unicyclic graphs}

In this section we determine the unicyclic graphs on $n$ vertices
with maximum reverse degree distances when girth, number of pendant
vertices and maximum degree are given respectively.

\begin{Lemma}
\label{th2} For $n\ge6$, $3\le m\le n-2$ and $2\le d<n-\left\lfloor
\frac{m+1}{2}\right\rfloor$,  $\r (U_{n,m,d}) < \r (U_{n,m,d+1})$.
\end{Lemma}

\begin{Proof}
Let $h=n-d-\lfloor\frac{m+1}{2}\rfloor$. Let $u_2$ be a pendant
vertex attached to $v_0$ different from $u$ in $U_{n,m,d}$ if $h\ge
2$. Recall that $U_{n,m,d}=U_{n,m,d}(\gamma,\theta)$. Note that we
may obtain $U_{n,m,d+1}(\gamma+1,\theta)$ from
$U_{n,m,d}(\gamma,\theta)-\{uv_0\}+\{uu_0\}$. Let
$G_1=U_{n,m,d+1}(\gamma+1,\theta)$ and
$G_2=U_{n,m,d}(\gamma,\theta)$. If $\theta\ge 1$ and $h\ge2$, then
$D_{G_1}\left(v_{\lfloor\frac{m}{2}\rfloor}\right)-D_{G_2}\left(v_{\lfloor\frac{m}{2}\rfloor}\right)=D_{G_1}(u_1)-D_{G_2}(u_1)$,
and thus
\begin{eqnarray*}
&&D'(U_{n,m,d+1}(\gamma+1,\theta))-D'(U_{n,m,d}(\gamma,\theta))\\
&=&4[W(G_1)-W(G_2)]-(h-1)[D_{G_1}(u_2)-D_{G_2}(u_2)]\\
&&-[D_{G_1}(u)-D_{G_2}(u)]+hD_{G_1}(v_0)-(h+1)D_{G_2}(v_0)+D_{G_2}(u_0)\\
&=&4\gamma(n-\gamma-2)-\gamma(h-1)-\gamma(n-\gamma-2)+\gamma(n-\gamma+h-1)\\
&=&-4\gamma^2+2(2n-3)\gamma.
\end{eqnarray*}
By similar calculation, the equality above holds also if $\theta=0$
or $h=1$. Thus
\begin{eqnarray*}
&&\r (U_{n,m,d+1}(\gamma+1,\theta))-\r (U_{n,m,d}(\gamma,\theta))\\
&=&2n(n-1)-[D'(U_{n,m,d+1}(\gamma+1,\theta))-D'(U_{n,m,d}(\gamma,\theta))]\\
&=&4\gamma^2-2(2n-3)\gamma+2n^2-2n\\
&\ge&4\left(\frac{2n-3}{4}\right)^2-2(2n-3)\cdot\frac{2n-3}{4}+2n^2-2n\\
&=&n^2+n-\frac{9}{4}>0.
\end{eqnarray*}
By Corollary \ref{Cor}, we have $\r
(U_{n,m,d+1}\left(\gamma+1,\theta)\right) \leq \r
\left(U_{n,m,d+1}\right)$, then the result follows.
\end{Proof}

\begin{Theorem}
Let $G$ be a  unicyclic graph with $n$ vertices and girth $m$, where
$n\ge 6$, $3\le m\le n-2$. Then $\r (G) \leq \r
\left(U_{n,m,n-\left\lfloor\frac{m+1}{2}\right\rfloor}\right)$ with
equality if and only if $G =
U_{n,m,n-\left\lfloor\frac{m+1}{2}\right\rfloor}$.
\end{Theorem}

\begin{Proof}
Let $d$ be the diameter of $G$.  Then $2\le d \leq
n-\lfloor\frac{m+1}{2}\rfloor$. If
$d=n-\lfloor\frac{m+1}{2}\rfloor$, then $\alpha=0$, and by Theorem
\ref{th1} (iv), the result follows. If $d=2$, then
$G=U_{n,3,2}(1,0)$, and by Lemma \ref{th2}, we have $\r
(U_{n,3,2}(1,0))<\r (U_{n,3,3})$.   If $3\le
d<n-\lfloor\frac{m+1}{2}\rfloor$, then by Corollary \ref{Cor} and
Lemma \ref{th2}, $\r (G) \leq \r (U_{n,m,d})<\r
\left(U_{n,m,n-\lfloor\frac{m+1}{2}\rfloor}\right)$.
\end{Proof}


\begin{Lemma}\cite{DZ}\label{l9}
For $n\ge5$, $3\le m\le n-1$ and $3\le d \le
n-\lfloor\frac{m+1}{2}\rfloor$, let
$h=n-d-\left\lfloor\frac{m+1}{2}\right\rfloor$. Then
\begin{eqnarray*}
&&W\left(U_{n,m,d}(a,b)\right)\\
&=&\left(a+b+\frac{m}{2}\right)\left\lfloor\frac{m^2}{4}\right\rfloor
+{a+1 \choose 3}+{b+1 \choose 3}\\
&&+m\left[{a+1 \choose 2}+{b+1 \choose 2}\right]
+\frac{1}{2}ab\left(2\left\lfloor\frac{m}{2}\right\rfloor+a+b+2\right)\\
&&+h\left[\left\lfloor\frac{m^2}{4}\right\rfloor+m
+\frac{1}{2}a(a+3)+\frac{1}{2}b\left(2\left\lfloor\frac{m}{2}\right\rfloor+b+3\right)\right]+h(h-1),
\end{eqnarray*}
where $a, b$ are integers with
$a+b=d-\left\lfloor\frac{m}{2}\right\rfloor$, $a\ge b\ge 0$ and
$a\ge 1$.
\end{Lemma}

By simple calculation, we have
\begin{Lemma}\label{l12}
For $G=U_{n,m,d}(a,b)$ with $n\ge5$, $3\le m\le n-1$ and $3\le d \le
n-\lfloor\frac{m+1}{2}\rfloor$, let
$h=n-d-\left\lfloor\frac{m+1}{2}\right\rfloor$. Then
\begin{eqnarray*}
D_G(v_0)&=&\left\lfloor\frac{m^2}{4}\right\rfloor+\frac{1}{2}a(a+1)
+\frac{1}{2}b\left(b+1+2\left\lfloor\frac{m}{2}\right\rfloor\right)+h,\\
D_G\left(v_{\left\lfloor\frac{m}{2}\right\rfloor}\right)&=&\left\lfloor\frac{m^2}{4}\right\rfloor+\frac{1}{2}a\left(a+1+2\left\lfloor\frac{m}{2}\right\rfloor\right)+\frac{1}{2}b(b+1)
+h\left(1+\left\lfloor\frac{m}{2}\right\rfloor\right),\\
D_G(u)&=&\left\lfloor\frac{m^2}{4}\right\rfloor+m+\frac{1}{2}a(a+3)
+\frac{1}{2}b\left(2\left\lfloor\frac{m}{2}\right\rfloor+b+3\right)+2(h-1),\\
D_G(u_0)&=&\left\lfloor\frac{m^2}{4}\right\rfloor+a\left[\frac{1}{2}(a-1)+m\right]+\frac{1}{2}b\left(2a+2\left\lfloor\frac{m}{2}\right\rfloor+b+1\right)+h(a+1),\\
D_G(u_1)&=&\left\lfloor\frac{m^2}{4}\right\rfloor+b\left[\frac{1}{2}(b-1)+m\right]
+\frac{1}{2}a\left(2b+2\left\lfloor\frac{m}{2}\right\rfloor+a+1\right)\\
&&+h\left(b+\left\lfloor\frac{m}{2}\right\rfloor+1\right).
\end{eqnarray*}
\end{Lemma}

\begin{Lemma}\label{l10}
Let $n$ and $m$ be integers with $5\le m\le n-1$. Let
$d=n-\lfloor\frac{m+1}{2}\rfloor$ and $a=n-m$. Then
\[
\r (U_{n,m,d}(a,0))<\r (U_{n,m-2,d+1}(a+2,0)).
\]
\end{Lemma}

\begin{Proof}
Let $G_1=U_{n,m-2,d+1}(a+2,0)$ and $G_2=U_{n,m,d}(a,0)$. Note that
$h=n-d-\left\lfloor\frac{m+1}{2}\right\rfloor=0$. By Lemmas \ref{l9}
and \ref{l12}, we have
\begin{eqnarray*}
&&D'(U_{n,m-2,d+1}(a+2,0))-D'(U_{n,m,d}(a,0))\\
&=&4[W(G_1)-W(G_2)]-[D_{G_1}(u_0)-D_{G_2}(u_0)]+[D_{G_1}(v_0)-D_{G_2}(v_0)]\\
&=&4\left[-\frac{3}{2}m^2+\left(n+\frac{9}{2}\right)m+\left\lfloor\frac{m^2}{4}\right\rfloor-2n-4\right]-(m-2)+(2n-3m+4)\\
&=&-6m^2+2(2n+7)m+4\left\lfloor\frac{m^2}{4}\right\rfloor-6n-10.
\end{eqnarray*}
Thus,
\begin{eqnarray*}
&&\r (U_{n,m-2,d+1}(a+2,0))- \r (U_{n,m,d}(a,0))\\
&=&2(n-1)n-[D'(U_{n,m-2,d+1}(a+2,0))-D'(U_{n,m,d}(a,0))]\\
&=&6m^2-2(2n+7)m-4\left\lfloor\frac{m^2}{4}\right\rfloor+2n^2+4n+10\\
&=&\left\{
\begin{array}{ll}
5m^2-2(2n+7)m+2n^2+4n+10 & \mbox{if $m$ is even}\\
5m^2-2(2n+7)m+2n^2+4n+11 & \mbox{if $m$ is odd}
\end{array}
\right.\\
&\ge&5m^2-2(2n+7)m+2n^2+4n+10\\
&\ge&5\cdot\left(\frac{2n+7}{5}\right)^2-2(2n+7)\cdot\frac{2n+7}{5}+2n^2+4n+10\\
&=&\frac{1}{5}(6n^2-8n+1)>0.
\end{eqnarray*}
Now the result follows.
\end{Proof}

\begin{Lemma}\label{l11}
Let $n$, $m$ and $d$ be integers with $5\le m\le n-2$, $3\le d \le
n-\lfloor\frac{m+1}{2}\rfloor$, and let $a$ and $b$ be integers with
$a+b=d-\left\lfloor\frac{m}{2}\right\rfloor$, $a\ge b\ge 0$ and
$a\ge 1$. Then
\[
\r \left(U_{n,m,d}(a,b)\right)<\r
\left(U_{n,m-2,d+1}(a+1,b+1)\right).
\]
\end{Lemma}

\begin{Proof}
Let $h=n-d-\lfloor\frac{m+1}{2}\rfloor$. Let
$G_1=U_{n,m-2,d+1}(a+1,b+1)$ and $G_2=U_{n,m,d}(a,b)$. If $b\ge 1$,
then by Lemmas \ref{l9} and \ref{l12}, we have
\begin{eqnarray*}
&&D'(U_{n,m-2,d+1}(a+1,b+1))-D'(U_{n,m,d}(a,b))\\
&=&4[W(G_1)-W(G_2)]\\
&&+(h+1)[D_{G_1}(v_0)-D_{G_2}(v_0)]+\left[D_{G_1}\left(v_{\lfloor\frac{m}{2}\rfloor-1}\right)-D_{G_2}\left(v_{\lfloor\frac{m}{2}\rfloor}\right)\right]\\
&&-h[D_{G_1}(u)-D_{G_2}(u)]-[D_{G_1}(u_0)-D_{G_2}(u_0)]-[D_{G_1}(u_1)-D_{G_2}(u_1)]\\
&=&4\left[h\left(2+a-\left\lfloor\frac{m+1}{2}\right\rfloor\right)-2+ab+
\left\lfloor\frac{m^2}{4}\right\rfloor\right.\\
&& \left. +\left\lfloor\frac{m}{2}\right\rfloor(a+b+1)+\frac{3}{2}m-\frac{m^2}{2}\right]\\
&&+(h+1)\left(2+a-\left\lfloor\frac{m+1}{2}\right\rfloor\right)+\left(2+b-\left\lfloor\frac{m+1}{2}\right\rfloor-h\right)\\
&&-h\left(2+a-\left\lfloor\frac{m+1}{2}\right\rfloor\right)-\left(b+\left\lfloor\frac{m}{2}\right\rfloor+h\right)-\left(a+\left\lfloor\frac{m}{2}\right\rfloor\right)\\
&=&-4b^2+4(n-m-2h)b+4n\left(\left\lfloor\frac{m}{2}\right\rfloor+h\right)-2m^2-8mh\\
&&-4(m-1)\left(\left\lfloor\frac{m}{2}\right\rfloor-1\right)+4\left\lfloor\frac{m^2}{4}\right\rfloor-4h^2+6h.\\
\end{eqnarray*}
If $b=0$, then $a=n-h-m$, and  by similar calculation, the equality
above holds also. Thus,
\begin{eqnarray*}
&&\r (U_{n,m-2,d+1}(a+1,b+1))-\r (U_{n,m,d}(a,b))\\
&=&2(n-1)n-[D'(U_{n,m-2,d+1}(a+1,b+1))-D'(U_{n,m,d}(a,b))]\\
&=&4b^2-4(n-m-2h)b\\
&&+2n^2-4n\left(\left\lfloor\frac{m}{2}\right\rfloor+h+\frac{1}{2}\right)+2m^2+8mh+4(m-1)\left(\left\lfloor\frac{m}{2}\right\rfloor-1\right)\\
&&-4\left\lfloor\frac{m^2}{4}\right\rfloor+4h^2-6h\\
&\ge&4\left(\frac{n-m-2h}{2}\right)^2-4(n-m-2h)\cdot\frac{n-m-2h}{2}\\
&&+2n^2-4n\left(\left\lfloor\frac{m}{2}\right\rfloor+h+\frac{1}{2}\right)+2m^2+8mh+4(m-1)\left(\left\lfloor\frac{m}{2}\right\rfloor-1\right)\\
&&-4\left\lfloor\frac{m^2}{4}\right\rfloor+4h^2-6h\\
&=&\left\{
\begin{array}{ll}
 2m^2+2(2h-3)m+n^2-2n+4-6h & \mbox{if $m$ is even}\\
 \begin{array}{l} 2m^2+2(2h-3)m+n^2-2n+4-6h\\+2n-2m+3\end{array} &
\mbox{if $m$ is odd}
\end{array}
\right.\\
&\ge&2m^2+2(2h-3)m+n^2-2n+4-6h\\
&\ge&2\cdot5^2+2\cdot5(2h-3)+n^2-2n+4-6h\\
&=&n^2-2n+24+14h>0.
\end{eqnarray*}
Now the result follows.
\end{Proof}

Let $\mathcal {U}(n,p)$ be the set of unicyclic graphs with $n$
vertices and $p$ pendant vertices, where $0 \leq p \leq n-3$. The
case $p=0$ is trivial.

Any graph in $\mathcal {U}(n,n-3)$ may be obtained by attaching
$n-3$ pendant vertices to vertices of a triangle, and then it is
easily seen that $U_{n,3,3}$ attains maximum reverse degree distance
in $\mathcal {U}(n,n-3)$.

\begin{Theorem}
\label{th3} Among graphs in $\mathcal {U}(n,p)$, where $n\ge6$ and
$1 \leq p \leq n-4$,
\begin{enumerate}

\item[(i)]
if $p=1$, then $U_{n,4,n-2}(n-4,0)$ is the unique graph with maximum
reverse degree distance;

\item[(ii)]
if $p=2$, then $U_{n,4,n-2}$ is the unique graph with maximum
reverse degree distance;

\item[(iii)]
if $p=3$ and $n=7$, then $U_{7,3,4}$ is the unique graph with
maximum reverse degree distance;

\item[(iv)]
if $p=3$ and $n>7$ is odd, then $U^k_{n,4,n-3}$ for $k=0,1$ are the
unique graphs with maximum reverse degree distance;

\item[(v)]
if $p=3$ and $n\ge6$ is even, or $4\le p\le n-4$, then $U_{n,4,n-p}$
is the unique graph with maximum reverse degree distance for
$\left\lfloor\frac{n-p-1}{2}\right\rfloor>\frac{n+4}{6}$,
$U_{n,3,n-p}$ and $U_{n,4,n-p}$ are the unique graphs with maximum
reverse degree distance for
$\left\lfloor\frac{n-p-1}{2}\right\rfloor=\frac{n+4}{6}$,
$U_{n,3,n-p}$ is the unique graph with maximum reverse degree
distance for
$\left\lfloor\frac{n-p-1}{2}\right\rfloor<\frac{n+4}{6}$.
\end{enumerate}
\end{Theorem}

\begin{Proof}
Obviously,  $\mathcal
{U}(n,1)=\{U_{n,m,n-\lfloor\frac{m+1}{2}\rfloor}(n-m,0):3\le m\le
n-1\}$. By Lemma \ref{l10}, $\r
(U_{n,m,n-\lfloor\frac{m+1}{2}\rfloor}(n-m,0))<\r
(U_{n,3,n-2}(n-3,0))$  for odd  $m>3$, and  $\r
(U_{n,m,n-\lfloor\frac{m+1}{2}\rfloor}(n-m,0))<\r
(U_{n,4,n-2}(n-4,0))$ for even $m>4$. Let $G_1=U_{n,4,n-2}(n-4,0)$
and $G_2=U_{n,3,n-2}(n-3,0)$. It is easily seen that
\begin{eqnarray*}
&&\r (U_{n,4,n-2}(n-4,0))-\r (U_{n,3,n-2}(n-3,0))\\
&=&D'(U_{n,3,n-2}(n-3,0))-D'(U_{n,4,n-2}(n-4,0))\\
&=&4[W(G_2)-W(G_1)]+[D_{G_2}(v_0)-D_{G_1}(v_0)]-[D_{G_2}(u_0)-D_{G_1}(u_0)]\\
&=&4(n-4)+(n-5)-1=5n-22>0.
\end{eqnarray*}
Then (i) follows.

Suppose that $2\le p\le n-4$. Let
 $G\in \mathcal {U}(n,p)$, and let $d$ and $m$ be respectively the diameter and girth of
$G$. A diametrical path contains at most
$\lfloor\frac{m}{2}\rfloor+1$ vertices on $C_m$ and two pendant
vertices, and thus at most
$(n-m-p)+\lfloor\frac{m}{2}\rfloor+1+2=n-p+3-\lfloor\frac{m+1}{2}\rfloor$
vertices in $G$. Thus $d\le n-p+2-\lfloor\frac{m+1}{2}\rfloor$.

By Corollary \ref{Cor} and Lemma
 \ref{th2}, $\r (G) \le
\r (U_{n,3,d})\le\r (U_{n,3,n-p})$ for $m=3$, and  $\r (G) \le \r
(U_{n,4,d})\le\r (U_{n,4,n-p})$ for $m=4$.

By Corollary \ref{Cor} and Lemmas \ref{th2} and \ref{l11}, if
$m\ge5$, then for $i=\lfloor\frac{m-3}{2}\rfloor$, we have
\begin{eqnarray*}
\r (G) &\le &
\r (U_{n,m,d})\le\r \left(U_{n,m,n-p+2-\lfloor\frac{m+1}{2}\rfloor}(\gamma,\theta)\right)\\
&<& \r \left(U_{n,m-2i,n-p+2-\lfloor\frac{m+1}{2}\rfloor
+i}\left(\gamma+i,\theta
+i\right)\right)\\
&\le &\r \left(U_{n,m-2i,n-p+2-\lfloor\frac{m+1}{2}\rfloor+i}\right)\\
&=&\r \left(U_{n,m-2i,n-p}\right).
\end{eqnarray*}
Thus, $\r (G)<\r (U_{n,3,n-p})$ for odd $m>3$, and $\r (G)<\r
(U_{n,4,n-p})$ for even $m>4$.

Note that
$U_{n,4,n-p}=U_{n,4,n-p}\left(\frac{n-p-2}{2},\frac{n-p-2}{2}\right)$
if $p=2,3$ and $n-p$ is even,
$U_{n,4,n-p}=U_{n,4,n-p}\left(\left\lfloor\frac{n-p}{2}\right\rfloor,\left\lfloor\frac{n-p-3}{2}\right\rfloor\right)$
if $p=2,3$ and $n-p$ is odd or $4\le p\le n-4$, and
$U_{n,3,n-p}=U_{n,3,n-p}\left(\left\lfloor\frac{n-p}{2}\right\rfloor,\left\lfloor\frac{n-p-1}{2}\right\rfloor\right)$.
Let $G_3=U_{n,4,n-p}$ and $G_4=U_{n,3,n-p}$.

Suppose that $p=2,3$ and $n-p$ is even. Note that
$D'(U^0_{n,4,n-p})=D'(U^1_{n,4,n-p})$ from Lemma \ref{l6}. It is
easily seen that
\begin{eqnarray*}
&&\r (U_{n,4,n-p})-\r (U_{n,3,n-p})\\
&=&D'(U_{n,3,n-p})-D'(U_{n,4,n-p})\\
&=&4[W(G_4)-W(G_3)]+(p-1)[D_{G_4}(v_0)-D_{G_3}(v_0)]\\
&&+\left[D_{G_4}\left(v_1\right)-D_{G_3}\left(v_2\right)\right]-(p-2)[D_{G_4}(u)-D_{G_3}(u)]\\
&&-[D_{G_4}(u_0)-D_{G_3}(u_0)]-[D_{G_4}(u_1)-D_{G_3}(u_1)]\\
&=&4\left(\frac{n-p-2}{2}-p+2\right)+(1-p)+(2-p)+(p-2)+(1-p)+(p-2)\\
&=&2n-7p+4=\left\{
\begin{array}{ll}
2n-10>0 & \mbox{if $p=2$ and $n\ge6$ is even}\\
2n-17>0 & \mbox{if $p=3$ and $n>7$ is odd}\\
-3      & \mbox{if $p=3$ and $n=7$},
\end{array}
\right.\\
\end{eqnarray*}
and thus (ii), (iii) and (iv) follow.

If  $p=n-4$, then $U_{n,4,n-p}=U_{n,4,4}(2,0)$,
$U_{n,3,n-p}=U_{n,3,4}(2,1)$, and it is easily checked that $\r
(U_{n,4,4})-\r (U_{n,3,4})=D'(U_{n,3,4})-D'(U_{n,4,4})=2-n<0$. If
$p=2,3$ and $n-p$ is odd, or $4\le p\le n-5$, then by similar
calculation as above, we have
\begin{eqnarray*}
\r (U_{n,4,n-p})-\r
(U_{n,3,n-p})=6\left\lfloor\frac{n-p-1}{2}\right\rfloor-n-4,
\end{eqnarray*}
and thus (v) follows easily.
\end{Proof}

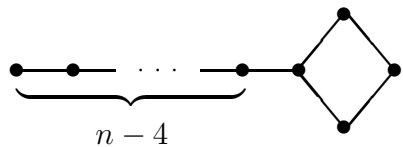
\begin{center}

\setlength{\unitlength}{0.75mm} \thicklines

{\begin{picture}(0,40)(-8,0)

\put(0,0){\begin{picture}(80,50)

\put(10,20){\circle*{2.5}}

 \put(1,30){\circle*{2.5}}
\put(1,10){\circle*{2.5}}\put(-7,20){\circle*{2.5}}

\put(10,20){\line(-5,6){8.5}} \put(10,20){\line(-5,-6){8.5}}
\put(1,30){\line(-5,-6){8.5}} \put(1,10){\line(-5,6){8.5}}

\put(-7,20){\line(-1,0){17.5}} \put(-39.5,20){\line(-1,0){17.5}}
\put(-17,20){\circle*{2.5}}
\put(-57,20){\circle*{2.5}}\put(-47,20){\circle*{2.5}}
\put(-29,20){\circle*{0.5}} \put(-32,20){\circle*{0.5}}
\put(-35,20){\circle*{0.5}}

\put(-43,7){${n-4}$}
\put(-57,17){$\underbrace{~~~~~~~~~~~~~~~~~~~~~~}$}

\put(-45,-5){$U_{n,4,n-2}(n-4,0)$}

\end{picture}}

\end{picture}}

\end{center}

\begin{center}

\setlength{\unitlength}{0.75mm} \thicklines

{\begin{picture}(50,40)(-8,0)

\put(0,0){\begin{picture}(80,50)

\put(10,20){\circle*{2.5}} \put(20,20){\circle*{2.5}}
\put(32,20){\circle*{0.5}} \put(35,20){\circle*{0.5}}
\put(38,20){\circle*{0.5}} \put(50,20){\circle*{2.5}}

\put(1,30){\circle*{2.5}} 
%
\put(1,10){\circle*{2.5}}\put(-7,20){\circle*{2.5}}

\put(10,20){\line(-5,6){8.5}} \put(10,20){\line(-5,-6){8.5}}
\put(1,30){\line(-5,-6){8.5}} \put(1,10){\line(-5,6){8.5}}

\put(10,20){\line(+1,0){17.5}} 
\put(42.5,20){\line(+1,0){17.5}} \put(60,20){\circle*{2.5}}


\put(-7,20){\line(-1,0){17.5}} 
\put(-39.5,20){\line(-1,0){17.5}} \put(-17,20){\circle*{2.5}}
\put(-57,20){\circle*{2.5}}\put(-47,20){\circle*{2.5}}
\put(-29,20){\circle*{0.5}} \put(-32,20){\circle*{0.5}}
\put(-35,20){\circle*{0.5}}

\put(-44,7){${\left\lceil\frac{n-4}{2}\right\rceil}$}
\put(-57,17){$\underbrace{~~~~~~~~~~~~~~~~~~~~~~}$}
\put(34,7){${\left\lfloor\frac{n-4}{2}\right\rfloor}$}
\put(20,17){$\underbrace{~~~~~~~~~~~~~~~~~~~~~~}$}

\put(-10,-5){$U_{n,4,n-2}$}

\end{picture}}
\end{picture}}

\end{center}

\begin{center}

\setlength{\unitlength}{0.75mm} \thicklines

{\begin{picture}(50,40)(-8,0)

\put(0,0){\begin{picture}(80,50)

\put(10,20){\circle*{2.5}} \put(20,20){\circle*{2.5}}
\put(32,20){\circle*{0.5}} \put(35,20){\circle*{0.5}}
\put(38,20){\circle*{0.5}} \put(50,20){\circle*{2.5}}

\put(1,30){\circle*{2.5}} 
%
\put(1,10){\circle*{2.5}}\put(-7,20){\circle*{2.5}}

\put(10,20){\line(-5,6){8.5}} \put(10,20){\line(-5,-6){8.5}}
\put(1,30){\line(-5,-6){8.5}} \put(1,10){\line(-5,6){8.5}}

\put(10,20){\line(+1,0){17.5}} 
\put(42.5,20){\line(+1,0){17.5}} \put(60,20){\circle*{2.5}}


\put(-7,20){\line(-1,0){17.5}} 
\put(-39.5,20){\line(-1,0){17.5}} \put(-17,20){\circle*{2.5}}
\put(-57,20){\circle*{2.5}}\put(-47,20){\circle*{2.5}}
\put(-29,20){\circle*{0.5}} \put(-32,20){\circle*{0.5}}
\put(-35,20){\circle*{0.5}}

\put(-7,20){\line(-3,4){8}}
\put(-15,30){\circle*{2.5}}

\put(-40,7){${\frac{n-5}{2}}$}
\put(-57,17){$\underbrace{~~~~~~~~~~~~~~~~~~~~~~}$}
\put(36,7){${\frac{n-5}{2}}$}
\put(20,17){$\underbrace{~~~~~~~~~~~~~~~~~~~~~~}$}

\put(-25,-5){$U_{n,4,n-3}$ (for odd $n$)}

\end{picture}}
\end{picture}}

\end{center}

\vskip 1cm

\begin{center}

\setlength{\unitlength}{0.75mm} \thicklines

{\begin{picture}(50,40)(-8,0)

\put(0,0){\begin{picture}(80,50)

\put(10,20){\circle*{2.5}} \put(20,20){\circle*{2.5}}
\put(32,20){\circle*{0.5}} \put(35,20){\circle*{0.5}}
\put(38,20){\circle*{0.5}} \put(50,20){\circle*{2.5}}

\put(1,30){\circle*{2.5}} 
%
\put(1,10){\circle*{2.5}}\put(-7,20){\circle*{2.5}}

\put(10,20){\line(-5,6){8.5}} \put(10,20){\line(-5,-6){8.5}}
\put(1,30){\line(-5,-6){8.5}} \put(1,10){\line(-5,6){8.5}}

\put(10,20){\line(+1,0){17.5}} 
\put(42.5,20){\line(+1,0){17.5}} \put(60,20){\circle*{2.5}}


\put(-7,20){\line(-1,0){17.5}} 
\put(-39.5,20){\line(-1,0){17.5}} \put(-17,20){\circle*{2.5}}
\put(-57,20){\circle*{2.5}}\put(-47,20){\circle*{2.5}}
\put(-29,20){\circle*{0.5}} \put(-32,20){\circle*{0.5}}
\put(-35,20){\circle*{0.5}}

\put(-7,20){\line(-3,4){8}}\put(-7,20){\line(1,5){2.1}}
\put(-15,30){\circle*{2.5}}\put(-5,30){\circle*{2.5}}
\put(-12,30){\circle*{0.5}}\put(-10,30){\circle*{0.5}}\put(-8,30){\circle*{0.5}}

\put(-16,37){${p-2}$} \put(-16,32){$\overbrace{~~~~~~}$}
\put(-43,7){${\left\lfloor\frac{n-p}{2}\right\rfloor}$}
\put(-57,17){$\underbrace{~~~~~~~~~~~~~~~~~~~~~~}$}
\put(30,7){${\left\lfloor\frac{n-p-3}{2}\right\rfloor}$}
\put(20,17){$\underbrace{~~~~~~~~~~~~~~~~~~~~~~}$}

\put(-55,-5){$U_{n,4,n-p}$ (for $p=3$ and even $n$, or $4\le p\le
n-4$)}

\end{picture}}
\end{picture}}

\end{center}

\vskip 1cm

\begin{center}

\setlength{\unitlength}{0.75mm} \thicklines

{\begin{picture}(50,40)(-8,0)

\put(0,0){\begin{picture}(80,50)

\put(10,20){\circle*{2.5}} \put(20,20){\circle*{2.5}}
\put(32,20){\circle*{0.5}} \put(35,20){\circle*{0.5}}
\put(38,20){\circle*{0.5}} \put(50,20){\circle*{2.5}}

\put(1,30){\circle*{2.5}} 
%
\put(-7,20){\circle*{2.5}}

\put(10,20){\line(-5,6){8.5}} 
\put(1,30){\line(-5,-6){8.5}} 
\put(10,20){\line(-1,0){16}}

\put(10,20){\line(+1,0){17.5}} 
\put(42.5,20){\line(+1,0){17.5}} \put(60,20){\circle*{2.5}}


\put(-7,20){\line(-1,0){17.5}} 
\put(-39.5,20){\line(-1,0){17.5}} \put(-17,20){\circle*{2.5}}
\put(-57,20){\circle*{2.5}}\put(-47,20){\circle*{2.5}}
\put(-29,20){\circle*{0.5}} \put(-32,20){\circle*{0.5}}
\put(-35,20){\circle*{0.5}}

\put(-7,20){\line(-3,4){8}}\put(-7,20){\line(1,5){2.1}}
\put(-15,30){\circle*{2.5}}\put(-5,30){\circle*{2.5}}
\put(-12,30){\circle*{0.5}}\put(-10,30){\circle*{0.5}}\put(-8,30){\circle*{0.5}}

\put(-16,37){${p-2}$} \put(-16,32){$\overbrace{~~~~~~}$}
\put(-43,7){${\left\lfloor\frac{n-p}{2}\right\rfloor}$}
\put(-57,17){$\underbrace{~~~~~~~~~~~~~~~~~~~~~~}$}
\put(30,7){${\left\lfloor\frac{n-p-1}{2}\right\rfloor}$}
\put(20,17){$\underbrace{~~~~~~~~~~~~~~~~~~~~~~}$}

\put(-10,-5){$U_{n,3,n-p}$}

\end{picture}}
\end{picture}}

\end{center}

\vspace{10mm}
\begin{center} Fig. 2.   The graphs in Theorem 3: $U_{n,4,n-2}(n-4,0)$, $U_{n,4,n-2}$,
$U_{n,4,n-3}$ (for odd $n$), $U_{n,4,n-p}$ (for $p=3$ and even $n$,
or $4\le p\le n-4$), and $U_{n,3,n-p}$\,.
\end{center}

\vspace{10mm}

Let $\mathfrak{U}(n,\Delta)$ be the set of unicyclic graphs with $n$
vertices and maximum degree $\Delta$, where $2 \leq \Delta \leq
n-1$. The cases $\Delta=2, n-1$ are trivial.

It is easily checked that $U_{n,3,3}$ attains maximum reverse degree
distance in $\mathfrak{U}(n, n-2)$.

\begin{Theorem} \label{th4}
Among graphs in $\mathfrak{U}(n,\Delta)$, where $n\ge6$ and $3\leq
\Delta \leq n-3$,

\begin{enumerate}

\item[(i)]
if $\Delta=3$, then $U_{n,4,n-2}$ is the unique graph with maximum
reverse degree distance;

\item[(ii)]
if $\Delta=4$ and $n=7$, then $U_{7,3,4}$ is the unique graph with
maximum reverse degree distance;

\item[(iii)]
if $\Delta=4$ and $n>7$ is odd, then $U^k_{n,4,n-3}$ for $k=0,1$ are
the unique graphs with maximum reverse degree distance;

\item[(iv)]
if $\Delta=4$ and $n\ge6$ is even, or $5\le \Delta\le n-3$, then
$U_{n,4,n-\Delta+1}$ is the unique graph with maximum reverse degree
distance for
$\left\lfloor\frac{n-\Delta}{2}\right\rfloor>\frac{n+4}{6}$,
$U_{n,3,n-\Delta+1}$ and $U_{n,4,n-\Delta+1}$ are the unique graphs
with maximum reverse degree distance for
$\left\lfloor\frac{n-\Delta}{2}\right\rfloor=\frac{n+4}{6}$,
$U_{n,3,n-\Delta+1}$ is the unique graph with maximum reverse degree
distance for
$\left\lfloor\frac{n-\Delta}{2}\right\rfloor<\frac{n+4}{6}$.

\end{enumerate}
\end{Theorem}

\begin{Proof}
Let $d$ be the diameter and let $u$ be a vertex of degree $\Delta$
in $G$.  A diametrical path $P(G)$ contains at most
$\lfloor\frac{m}{2}\rfloor+1$ vertices of the cycle $C_m$ and two
neighbors of $u$, but $P(G)$ can not contain these vertices at the
same time. Note that the cycle $C_m$ contains at most two neighbors
of $u$. Thus $d+1<n-(m+\Delta-2)+\lfloor\frac{m}{2}\rfloor+1+2$ and
then $d\le n-\Delta+3-\lfloor\frac{m+1}{2}\rfloor$.

By Corollary \ref{Cor} and Lemmas \ref{th2} and \ref{l11}, we have:
$\r (G) \le \r (U_{n,3,d})$ $\le\r\left(U_{n,3,n-\Delta+1}\right)$
for $m=3$, $\r (G) \le \r \left(U_{n,4,d})\le\r
(U_{n,4,n-\Delta+1}\right)$ for $m=4$, and
\begin{eqnarray*}
\r (G) &\le &
\r (U_{n,m,d})\le\r \left(U_{n,m,n-\Delta+3-\lfloor\frac{m+1}{2}\rfloor}(\gamma,\theta)\right)\\
&<&\r \left(U_{n,m-2i,n-\Delta+3-\lfloor\frac{m+1}{2}\rfloor+i}\left(\gamma+i,\theta+i\right)\right)\\
&\le
&\r \left(U_{n,m-2i,n-\Delta+3-\lfloor\frac{m+1}{2}\rfloor+i}\right)\\
&=& \r \left(U_{n,m-2i,n-\Delta+1}\right)
\end{eqnarray*}
for $m\ge 5$, where $i=\lfloor\frac{m-3}{2}\rfloor$.

Thus, $\r (G)<\r\left(U_{n,3,n-\Delta+1}\right)$ for odd $m>3$ and
 $\r (G)<\r\left(U_{n,4,n-\Delta+1}\right)$ for even $m>4$.
Now the theorem follows by similar arguments as in the proof of
Theorem~\ref{th3}.
\end{Proof}

Finally, we give the values of the maximum reverse degree distances
in Theorem \ref{th3} and \ref{th4}.

(i) For $U_{n,4,n-2}(n-4,0)$,
\begin{eqnarray*}
D'(U_{n,4,n-2}(n-4,0))&=&4W(U_{n,4,n-2}(n-4,0))+(3-2)D_{U_{n,4,n-2}(n-4,0)}(v_0)\\
&&+(1-2)D_{U_{n,4,n-2}(n-4,0)}(u_0)\\
&=&4W(U_{n,4,n-2}(n-4,0))-3(n-4)\\
&=&\frac{2}{3}n^3-\frac{35}{3}n+36,
\end{eqnarray*}
and thus,
\begin{eqnarray*}
^rD'(U_{n,4,n-2}(n-4,0))&=&2(n-1)n(n-2)-D'(U_{n,4,n-2}(n-4,0))\\
&=&\frac{4}{3}n^3-6n^2+\frac{47}{3}n-36.
\end{eqnarray*}

(ii) For $U_{n,4,n-2}$,
\begin{eqnarray*}
D'(U_{n,4,n-2})&=&4W(U_{n,4,n-2})+(3-2)D_{U_{n,4,n-2}}(v_0)+(1-2)D_{U_{n,4,n-2}}(u_0)\\
&&+(3-2)D_{U_{n,4,n-2}}(v_{\lfloor\frac{m}{2}\rfloor})+(1-2)D_{U_{n,4,n-2}}(u_1)\\
&=&4W(U_{n,4,n-2}(n-4,0))-3(n-4)\\
&=&\left\{
\begin{array}{ll}
\frac{2}{3}n^3-\frac{3}{2}n^2+\frac{1}{3}n+12 & \mbox{if $n$ is even}\\
\frac{2}{3}n^3-\frac{3}{2}n^2+\frac{1}{3}n+\frac{27}{2} & \mbox{if
$n$ is odd},
\end{array}
\right.
\end{eqnarray*}
and thus
\begin{eqnarray*}
^rD'(U_{n,4,n-2})&=&2(n-1)n(n-2)-D'(U_{n,4,n-2})\\
&=&\left\{
\begin{array}{ll}
\frac{4}{3}n^3-\frac{9}{2}n^2+\frac{11}{3}n-12 & \mbox{if $n$ is even}\\
\frac{4}{3}n^3-\frac{9}{2}n^2+\frac{11}{3}n-\frac{27}{2} & \mbox{if
$n$ is odd}.
\end{array}
\right.
\end{eqnarray*}

(iii) For $U_{n,4,n-3}$ with odd $n$,
\begin{eqnarray*}
D'(U_{n,4,n-3})&=&4W(U_{n,4,n-3})+(4-2)D_{U_{n,4,n-3}}(v_0)+(1-2)D_{U_{n,4,n-3}}(u_0)\\
&&+(1-2)D_{U_{n,4,n-3}}(u)+(3-2)D_{U_{n,4,n-3}}(v_{\lfloor\frac{m}{2}\rfloor})\\
&&+(1-2)D_{U_{n,4,n-3}}(u_1)\\
&=&4W(U_{n,4,n-3})-\frac{1}{2}n^2+\frac{19}{2}\\
&=&\frac{2}{3}n^3-\frac{5}{2}n^2+\frac{10}{3}n+\frac{47}{2}
\end{eqnarray*}
and thus
\begin{eqnarray*}
^rD'(U_{n,4,n-3})&=&2(n-1)n(n-3)-D'(U_{n,4,n-3})\\
&=&\frac{4}{3}n^3-\frac{11}{2}n^2+\frac{8}{3}n-\frac{47}{2}.
\end{eqnarray*}

(iv) For $U_{n,4,n-p}$,
\begin{eqnarray*}
&&D'(U_{n,4,n-p})\\
&=&4W(U_{n,4,n-p})+[(p+1)-2]D_{U_{n,4,n-p}}(v_0)\\
&&+(1-2)D_{U_{n,4,n-p}}(u_0)+(p-2)(1-2)D_{U_{n,4,n-p}}(u)\\
&&+(3-2)D_{U_{n,4,n-p}}(v_{\lfloor\frac{m}{2}\rfloor})+(1-2)D_{U_{n,4,n-p}}(u_1)\\
&=&\left\{
\begin{array}{ll}
\frac{2}{3}n^3-\left(p-\frac{1}{2}\right)n^2+\left(3p-\frac{17}{3}\right)n+\frac{p^3}{3}
-\frac{p^2}{2}+\frac{11p}{3}+10 & \mbox{if $n-p$ is even}\\
\frac{2}{3}n^3-\left(p-\frac{1}{2}\right)n^2+\left(3p-\frac{17}{3}\right)n
+\frac{p^3}{3}-\frac{p^2}{2}+\frac{14p}{3}+\frac{7}{2}
& \mbox{if $n-p$ is odd},
\end{array}
\right.
\end{eqnarray*}
and thus
\begin{eqnarray*}
&& ^rD'(U_{n,4,n-p})\\
&=&2(n-1)n(n-p)-D'(U_{n,4,n-p})\\
&=&\left\{
\begin{array}{ll}
\frac{4}{3}n^3-\left(p+\frac{5}{2}\right)n^2-\left(p-\frac{17}{3}\right)n-\frac{p^3}{3}+\frac{p^2}{2}-\frac{11p}{3}-10 & \mbox{if $n$ is even}\\
\frac{4}{3}n^3-\left(p+\frac{5}{2}\right)n^2-\left(p-\frac{17}{3}\right)n-\frac{p^3}{3}+\frac{p^2}{2}-\frac{14p}{3}-\frac{7}{2}
& \mbox{if $n$ is odd}.
\end{array}
\right.
\end{eqnarray*}

(v) For $U_{n,3,n-p}$,
\begin{eqnarray*}
&&D'(U_{n,3,n-p})\\
&=&4W(U_{n,3,n-p})+[(p+1)-2]D_{U_{n,3,n-p}}(v_0)+(1-2)D_{U_{n,3,n-p}}(u_0)\\
&&+(p-2)(1-2)D_{U_{n,3,n-p}}(u)+(3-2)D_{U_{n,3,n-p}}(v_{\lfloor\frac{m}{2}\rfloor})+(1-2)D_{U_{n,3,n-p}}(u_1)\\
&=&\left\{
\begin{array}{ll}
\frac{2}{3}n^3-\left(p-\frac{1}{2}\right)n^2+\left(3p-\frac{11}{3}\right)n+\frac{p^3}{3}-\frac{p^2}{2}+\frac{2p}{3} & \mbox{if $n-p$ is even}\\
\frac{2}{3}n^3-\left(p-\frac{1}{2}\right)n^2+\left(3p-\frac{11}{3}\right)n+\frac{p^3}{3}-\frac{p^2}{2}+\frac{5p}{3}-\frac{7}{2}
& \mbox{if $n-p$ is odd},
\end{array}
\right.
\end{eqnarray*}
and thus
\begin{eqnarray*}
&&^rD'(U_{n,3,n-p})\\
&=&2(n-1)n(n-p)-D'(U_{n,3,n-p})\\
&=&\left\{
\begin{array}{ll}
\frac{4}{3}n^3-\left(p+\frac{5}{2}\right)n^2-\left(p-\frac{11}{3}\right)n-\frac{p^3}{3}+\frac{p^2}{2}-\frac{2p}{3} & \mbox{if $n-p$ is even}\\
\frac{4}{3}n^3-\left(p+\frac{5}{2}\right)n^2-\left(p-\frac{11}{3}\right)n-\frac{p^3}{3}+\frac{p^2}{2}-\frac{5p}{3}+\frac{7}{2}
& \mbox{if $n-p$ is odd}.
\end{array}
\right.
\end{eqnarray*}

\vspace{6mm}

\noindent {\bf Acknowledgements\/} This work was supported by the
National Natural Science Foundation of China
(no.~11071089).

\end{document}